%% file: conv_orient_europ_j_comb.tex
\documentclass [12pt]{elsart}
\usepackage{amsfonts, amsmath}
\usepackage{amssymb}
\usepackage{epsfig, subfigure, verbatim}

\renewcommand{\emph}[1]{{\it #1}}
\newcounter{segcount}

\newenvironment{segment}
{\refstepcounter{segcount}\vspace{5mm}
\noindent{\bf \thesegcount. }}
{}

\newenvironment{statementnumbered}[3][]
{\refstepcounter{segcount}\vspace{5mm}
\noindent{\bf\thesegcount. #2}#1{\bf.}\ {\sl #3}}
{\nolinebreak[4] \nopagebreak[4] $\hfill \square$}
\newenvironment{statement}[3][]
{\vspace{5mm}\noindent{\bf #2}#1{\bf.} {\sl #3}}
{\nolinebreak[4] \nopagebreak[4] $\hfill \square$}
\newenvironment{statementnoboxnumbered}[3][]
{\refstepcounter{segcount}\vspace{5mm}
\noindent{\bf\thesegcount. #2}#1{\bf.} {\sl #3}}
{}
\newenvironment{statementnobox}[3][]
{\vspace{5mm}\noindent{\bf #2}#1{\bf.} {\sl #3}}
{}
\newenvironment{definitionnumbered}[2][]
{\refstepcounter{segcount}\vspace{5mm}
\noindent{\bf\thesegcount. #2}#1{\bf.}}
{}

\newenvironment{definition}[2][]
{\vspace{5mm}\noindent{\bf #2}#1{\bf.}}
{}

\newenvironment{resultnumbered}[3][]
{\refstepcounter{segcount}\vspace{5mm}
\noindent{\bf\thesegcount. #2}#1{\bf.} {\sl #3}
\vskip5mm\noindent {\bf Proof: }}
{\nolinebreak[4] \nopagebreak[4] $\hfill \square$}

\newenvironment{result}[3][]
{\vspace{5mm}
\noindent{\bf#2}#1{\bf.} {\sl #3}
{\\ \bf Proof: }}
{\nolinebreak[4] \nopagebreak[4] $\hfill \square$}

\newenvironment{risultnumbered}[3][]
{\refstepcounter{segcount}\vspace{5mm}
\noindent{\bf\thesegcount. #2}#1{\bf.} {\sl #3}
{\nopagebreak[4] \noindent \bf Proof: }}
{$\hfill \square$}

\newenvironment{risult}[3][]
{\vspace{5mm}
\noindent{\bf#2}#1{\bf.} {\sl #3}
{\noindent \bf Proof: }}
{\nolinebreak[4] \nopagebreak[4] $\hfill \square$}

\newenvironment{risultnoboxnumbered}[3][]
{\refstepcounter{segcount}\vspace{5mm}
\noindent{\bf\thesegcount. #2}#1{\bf.} {\sl #3}
{\nopagebreak[4] \noindent \bf Proof: }}
{}

\newenvironment{risultnobox}[3][]
{\vspace{5mm}
\noindent{\bf#2}#1{\bf.} {\sl #3}
{\nopagebreak[4] \noindent \bf Proof: }}
{}

\newenvironment{resultnoboxnumbered}[3][]
{\refstepcounter{segcount}\vspace{5mm}
\noindent{\bf\thesegcount. #2}#1{\bf.} {\sl #3}
{\\ \bf Proof: }}
{}

\newenvironment{resultnobox}[3][]
{\vspace{5mm}\noindent{\bf #2}#1{\bf.} {\sl #3}
{\\ \bf Proof: }}
{}

\newcommand{\seg}{\begin{segment}}
\newcommand{\segend}{\end{segment}}

\newcommand{\stmtnum}{\begin{statementnumbered}}
\newcommand{\stmtnumend}{\end{statementnumbered}}

\newcommand{\stmt}{\begin{statement}}
\newcommand{\stmtend}{\end{statement}}

\newcommand{\stmtnoboxnum}{\begin{statementnoboxnumbered}}
\newcommand{\stmtnoboxnumend}{\end{statementnoboxnumbered}}

\newcommand{\stmtnobox}{\begin{statementnobox}}
\newcommand{\stmtnoboxend}{\end{statementnobox}}

\newcommand{\defnnum}{\begin{definitionnumbered}}
\newcommand{\defnnumend}{\end{definitionnumbered}}

\newcommand{\defnn}{\begin{definition}}
\newcommand{\defnnend}{\end{definition}}

\newcommand{\resnum}{\begin{resultnumbered}}
\newcommand{\resnumend}{\end{resultnumbered}}

\newcommand{\res}{\begin{result}}
\newcommand{\resend}{\end{result}}

\newcommand{\risnum}{\begin{risultnumbered}}
\newcommand{\risnumend}{\end{risultnumbered}}

\newcommand{\ris}{\begin{risult}}
\newcommand{\risend}{\end{risult}}

\newcommand{\risnoboxnum}{\begin{risultnoboxnumbered}}
\newcommand{\risnoboxnumend}{\end{risultnoboxnumbered}}

\newcommand{\risnobox}{\begin{risultnobox}}
\newcommand{\risnoboxend}{\end{risultnobox}}

\newcommand{\resnoboxnum}{\begin{resultnoboxnumbered}}
\newcommand{\resnoboxnumend}{\end{resultnoboxnumbered}}

\newcommand{\resnobox}{\begin{resultnobox}}
\newcommand{\resnoboxend}{\end{resultnobox}}

\begin{document}

\markboth{Orientable convexity number}{Orientable convexity number}

\begin{frontmatter}
\title{Orientable convexity, geodetic and hull numbers in graphs}

\author{Alastair Farrugia\thanksref{ccs}}
\thanks[ccs]{The author's studies are sponsored by the Canadian government, through a Canadian Commonwealth Scholarship.}
\ead{afarrugia@math.uwaterloo.ca}
\ead[url]{members.lycos.co.uk/afarrugia}
\address{
Dept. of Combinatorics \& Optimization
\\ University of Waterloo, Ontario, Canada, N2L 3G1
}

%\date{March, 2003}

\begin{abstract}
We prove three results conjectured or stated by Chartrand, Fink and Zhang [European J. Combin {\bf 21} (2000) 181--189, Disc.\ Appl.\ Math.\ {\bf 116} (2002) 115--126, and pre-print of ``The hull number of an oriented graph'']. For a digraph $D$, Chartrand et al.\ defined the geodetic, hull and convexity number --- $g(D)$, $h(D)$ and $con(D)$, respectively. For an undirected graph $G$, $g^{-}(G)$ and $g^{+}(G)$ are the minimum and maximum geodetic numbers over all orientations of $G$, and similarly for $h^{-}(G)$, $h^{+}(G)$, $con^{-}(G)$ and $con^{+}(G)$. Chartrand and Zhang gave a proof that $g^{-}(G) < g^{+}(G)$ for any connected graph with at least three vertices. We plug a gap in their proof, allowing us also to establish their conjecture that $h^{-}(G) < h^{+}(G)$.

If $v$ is an end-vertex, then in any orientation of $G$, $v$ is either a source or a sink. It is easy to see that graphs without end-vertices can be oriented to have no source or sink; we show that, in fact, we can avoid all extreme vertices. This proves another conjecture of Chartrand et al., that $con^{-}(G) < con^{+}(G)$ iff $G$ has no end-vertices. 
\end{abstract}

\begin{keyword}
graph \sep digraph \sep oriented graph \sep convex \sep geodesic \sep convexity number \sep hull number \sep geodetic number \sep transitively orientable
\end{keyword}

\end{frontmatter}

The aim of this paper is to establish the following results, for every connected graph $G$ with at least three vertices:

\begin{eqnarray}
g^{-}(G) &<& g^{+}(G) \label{g-eqn} \\
h^{-}(G) &<& h^{+}(G) \label{h-eqn} \\
con^{-}(G) &<& con^{+}(G) \textrm{ iff $G$ has no end-vertices.} \label{con-eqn}
\end{eqnarray}

Results (\ref{h-eqn}) and (\ref{con-eqn}) were conjectured by Chartrand, Fink and Zhang in~\cite{cfz03} and~\cite{cfz02}, respectively.
The first result was stated by Chartrand and Zhang in~\cite[Thm. 2.5]{cz00}, but there was a gap in their proof. They independently noticed this gap, and an alternative proof was found, but the correction we present in Section~\ref{sec-g-h} allows us to prove (\ref{g-eqn}) and (\ref{h-eqn}) simultaneously. We prove~(\ref{con-eqn}) in Section~\ref{sec-con}.

\section{Preliminaries}
Let $D = (V,A)$ be a digraph, and let $u$ and $v$ be vertices. A $u-v$ {\em geodesic} is a dipath from $u$ to $v$ with the least possible number of arcs. The {\em closed interval} $I[u,v]$ consists of $u$, $v$, and every vertex that is on some $u-v$ geodesic or on some $v-u$ geodesic (note that there may be no directed path at all from $u$ to $v$, or from $v$ to $u$). For a set $S \subseteq V(D)$, we define 
$I[S] := \bigcup_{u,v \in S} I[u,v]$,
and, for $k > 0$, $I^k[S] := I[I^{k-1}(S)]$, where $I^0[S] := S$.

A set $S$ is {\em convex} if $S = I[S]$, that is, every geodesic between every two vertices of $S$ lies in $S$. The {\em convex hull} $[S]$ of $S$ is the smallest convex set containing $S$; this is the intersection of all convex sets containing $S$, and also the limit of the sequence $S \subseteq I[S] \subseteq I^2[S] \subseteq \cdots.$

A {\em hull-set} of $D$ is a set $S \subseteq V$ for which $[S] = V$. If, moreover, $I[S] = V$, then $S$ is a {\em geodetic set}. The {\em hull number} of  $D$ is 

\[ h(D):= \min \{|S| \mid S \textrm{ is a hull-set of } D\},\]

while the {\em geodetic number} of  $D$ is 

\[ g(D):= \min \{|S| \mid S \textrm{ is a geodetic set of } D\}.\]

For an undirected graph $G$, an orientation $\overrightarrow{G}$ is a digraph obtained by giving each edge one of its two possible directions. The {\em lower} and {\em upper orientable hull numbers} are, respectively,

\begin{eqnarray*}
h^{-}(G) &:=& \min \{h(\overrightarrow{G}) \mid \overrightarrow{G} \textrm{ is an orientation of } G\}, \textrm{ and} \\
h^{+}(G) &:=& \max \{h(\overrightarrow{G}) \mid \overrightarrow{G} \textrm{ is an orientation of } G\}.
\end{eqnarray*}

The {\em lower} and {\em upper orientable geodetic numbers} $g^{-}(G)$ and $g^{+}(G)$ are defined similarly. 
\newline

Let $v$ be a vertex in a digraph $D = (V,A)$. Its {\em in-} and {\em out-neighbourhood} are $N^{-}(v) := \{u \mid uv \in A\}$  and $N^{+}(v) := \{w \mid vw \in A\}$, respectively. Its {\em in-} and {\em out-degree} are $id(v) := |N^{-}(v)|$ and $od(v) := |N^{+}(v)|$, respectively. 
If, for every $u \in N^{-}(v)$ and every $w \in N^{+}(v)$, $\overrightarrow{vw} \in A$, then $v$ is {\em extreme}. It is a {\em source} if $N^{-}(v) = \emptyset$, and a {\em sink} if $N^{+}(v) = \emptyset$. 
%An extreme vertex that is not a source or a sink is a {\em transitive} vertex. 

A graph that can be oriented so that every vertex is extreme is a {\em comparability} or {\em transitively orientable} graph.
A result that we will use repeatedly is the following, due to Chartrand et al.\ \cite[Prop.\ 2.1]{cfz02}, \cite[Prop.\ 1.3]{cfz03}:

\stmtnum{Proposition\label{extreme}}
{A vertex $v$ is extreme iff, for every $u$ and $w$ in $V$, $v$ is not an interior 
vertex of any $u-w$ geodesic. Therefore, $v$ is extreme iff $V -v$ is a convex set, iff 
$v$ is contained in every hull-set and every geodetic-set.
}
\stmtnumend
%\newline

\section{Orientable convexity numbers\label{sec-con}}
If $D = (V,A)$ is a digraph, the {\em convexity number} $con(D)$ is the size of the largest convex set $C \subsetneqq V$ ($V$ itself is always convex). For an undirected graph $G$, $con^{-}(G)$ and $con^{+}(G)$ are the minimum and maximum convexity numbers over all orientations of $G$. We are interested in whether $con^{-}(G) < con^{+}(G)$.

By Proposition~\ref{extreme}, if $D$ has an extreme vertex, then $con(D) = n-1$, where $n$ is the number of vertices. For any graph $G$, we can make an arbitrary vertex $v$ extreme by orienting all incident edges away from $v$, so we always have $con^{+}(G) = n-1$. Moreover, if $G$ contains an end-vertex $x$, then in every orientation $x$ is either a source or a sink; so in this case, $con^{-}(G) = n-1$~too.

If $G$ has no end-vertices, it is straightforward to find an orientation with no sources or sinks; the reader is encouraged to do so, and then try to generalise this to avoid all extreme vertices. We present a solution below.

%One solution involves ear decompositions for the $2$-connected case (the 
%$1$-connected case then following quickly), where each ear is taken to be as short as 
%possible, but there is a simpler proof that we present below.

Let some of the edges of $G$ be oriented. A vertex incident to some oriented edge is an {\em or-vertex}, short for {\em oriented vertex}. Note that a vertex $v$ is non-extreme iff there are arcs $\overrightarrow{uv}$ and $\overrightarrow{vw}$, such that $uw$ is either not present, or it is already oriented as $\overleftarrow{uw}$. No matter how the remaining undirected edges are oriented, $v$ remains non-extreme. 

\resnum{Theorem}
{A graph with minimum degree $2$ can be oriented so that all its vertices are 
non-extreme. Thus, for a connected graph $G$ with at least $3$ vertices, 
$con^{-}(G) < con^{+}(G)$ iff $G$ has no end-vertices.
}
Since $G$ has minimum degree $2$, it contains a cycle. Find a maximal set of 
%CHANGED 
edge-disjoint 
%CHANGED 
\emph{chordless} cycles, and orient their edges to make them directed cycles. We claim that every or-vertex $v$ is now non-extreme. If $v$ is on a triangle $uvw$, then $\overrightarrow{uv}$, $\overrightarrow{vw}$ and  $\overrightarrow{wu}$ are all arcs. Otherwise, $v$ is on a chordless cycle of length at least $4$, with neighbours, say, $u$ and $w$, where $uw \notin E(G)$. 

We now show that, if there are unoriented vertices, we can orient one or more while maintaining the property that all or-vertices are non-extreme.

Any unoriented vertex $u$ must be on a path $u_0, \ldots, u_{r+1}$ joining 
%CHANGED 
distinct 
%CHANGED 
or-vertices $u_0$ and $u_{r+1}$ (because the graph has minimum degree at least $2$, and our initial set of 
%CHANGED 
edge-disjoint 
%CHANGED 
cycles was chosen to be maximal).
Taking $r$ to be as small as possible ensures that the internal vertices $u_1, \ldots, u_r$ are all unoriented. Directing the path as $\overrightarrow{u_0 u_1}$, \ldots, $\overrightarrow{u_r u_{r-1}}$ ensures that $u_1, \ldots, u_{r}$ all have positive in- and out-degree. Moreover, if $r > 1$, then, for $1 \leq i \leq r$, 
$u_{i-1}u_{i+1} \notin E(G)$, and thus $u_i$ is non-extreme. 

If $r=1$, then we might have to orient differently as $u_0 u_2$ could be an edge of $G$.  If this edge is not oriented, we can orient it arbitrarily, since $u_0$ and $u_2$ are assumed to be already non-extreme. Without loss of generality, let it be oriented as $\overrightarrow{u_0 u_2}$; now orienting $u_0 u_1$ and $u_1 u_2$ as $\overleftarrow{u_0 u_1}$ and $\overleftarrow{u_1 u_2}$, ensures that $u_1$ is on a directed triangle and is thus non-extreme. 
\resnumend
%\newline

\section{Orientable geodetic and hull numbers\label{sec-g-h}}
Chartrand and Zhang's proof of~(\ref{g-eqn})
%\cite[Thm. 2.5]{cz00} 
essentially found a vertex $v_1$, and orientations $D_1$ and $D_2$ of $G$, such that if $S$ is a hull-set in $D_2$, then $I_{D_2}(S) \subseteq I_{D_1}(S-v_1)$ (this is Claim 1 in our own proof). Moreover, $v_1$ was a source in $D_2$, and was thus contained in every hull-set. By taking $S$ to be a minimum geodetic set for $D_2$, we immediately get $g^{-}(G) < g^{+}(G)$. With slightly more work (Claim 2 in our proof), we also get $h^{-}(G) < h^{+}(G)$, proving Conjecture 3.10 of~\cite{cfz03}.

Chartrand and Zhang  stated their result only for orientable geodetic numbers, as they did not include Claim 2. Moreover, they oriented $G[U]$ arbitrarily (where $U$ is defined in the proof). The path of length four (for example) shows that this does not always work, and their alternative proof did not extend to showing $h^{-}(G) < h^{+}(G)$. There is, however, an orientation of $G[U]$ that will rescue the original proof, as we show below.

\resnum{Theorem\label{g-h-eqn}}
{For any connected graph $G$ with at least three vertices, $g^{-}(G) < g^{+}(G)$ and $h^{-}(G) < h^{+}(G)$.
}
If $G$ is a complete graph with vertices $v_1, \ldots, v_n$, we first orient $G$ transitively (that is, $v_i \rightarrow v_j$ iff $i < j$). Since every vertex is extreme, this orientation shows that $g^{+}(G) = n = h^{+}(G)$. Reversing the orientation of $v_1 v_2$, \ldots, $v_{n-1} v_n$ makes $\{v_1, v_2\}$ a geodetic set; thus $g^{-}(G) = 2 = h^{-}(G)$. 

\begin{center}
\begin{figure}\label{fig-geo-hull}
\input{geo-hull.pstex_t}
\caption{The orientations $D_1$ and $D_2$ of $G$.}
\end{figure}
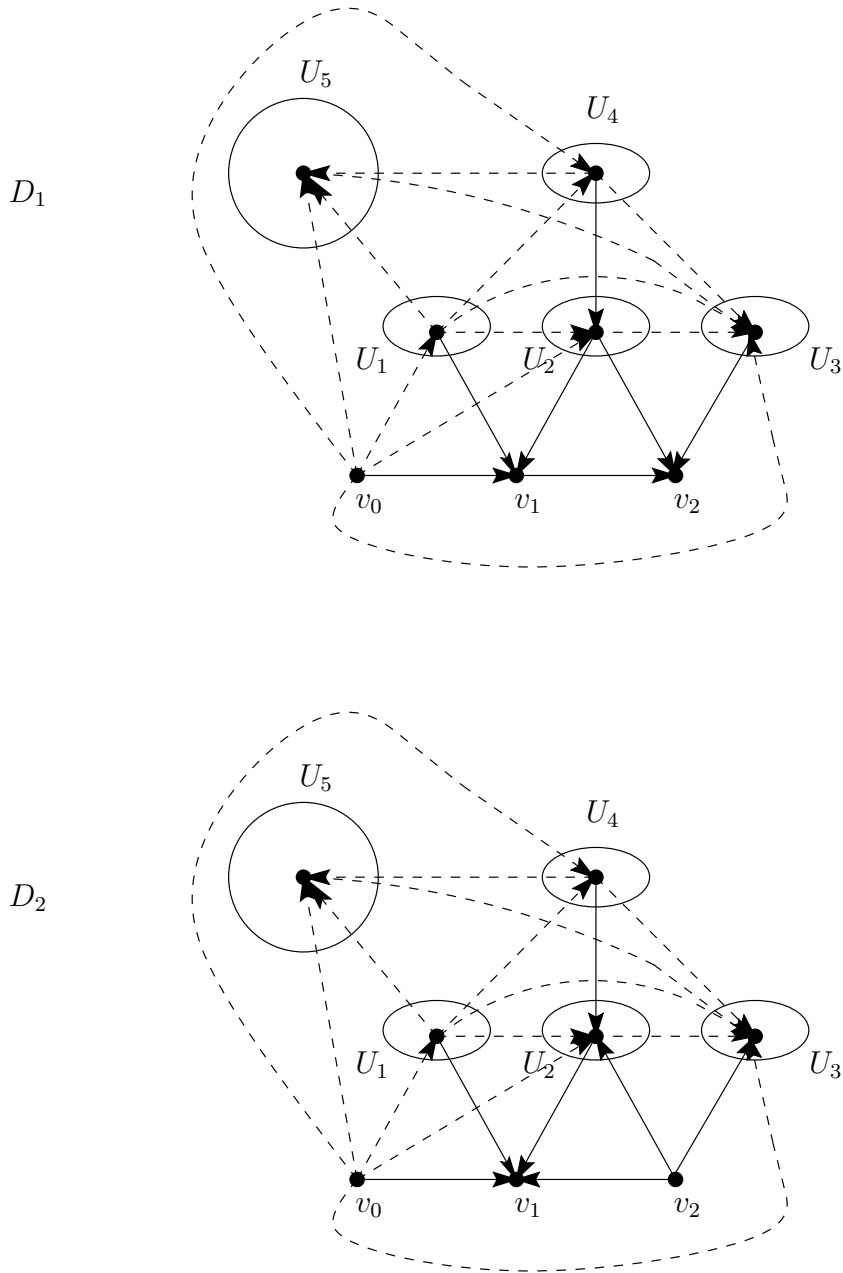
\end{center}
\vspace{-24pt}

If $G$ is not complete, then we can find vertices $v_0, v_1, v_2$ that induce a path of length two. Figure~\ref{fig-geo-hull} shows all the adjacencies (solid lines) and possible adjacencies (dashed lines) in $G$, where the $U_i$'s are defined as follows. 
%Let $G = (V,E)$ be a graph. 
For a set $C \subseteq V(G)$, $N(C)$ is the set $\{v \in V \mid \exists c \in C, vc \in E\}$.
\begin{eqnarray*}
U &:=& V(G) \setminus \{v_0, v_1, v_2\}, \\
U_1 &:=& U \cap (N(v_1) \setminus N(v_2)), \\
U_2 &:=& U \cap (N(v_1) \cap N(v_2)), \\
U_3 &:=& U \cap (N(v_2) \setminus N(v_1)), \\
U_4 &:=& (U \cap N(U_2)) \setminus (U_1 \cup U_2 \cup U_3), \textrm{ and} \\
U_5 &:=& U \setminus (U_1 \cup U_2 \cup U_3 \cup U_4).
\end{eqnarray*}

Let $D_2$ be the digraph\footnote
{The labeling is chosen to be consistent with Chartrand and Zhang, but I prefer to describe $D_2$ before $D_1$.} 
obtained by orienting $G$ as follows. We orient an edge $xy$ from $x$ to $y$ if one of the following conditions holds:
%\newline

\begin{center}
$\begin{array}{lll}
\multicolumn{2}{l}{x \in \{v_0, v_2\},} \\
&& y = v_1, \\
%x \in U_1 &\textrm{and}& y \in U_2 \cup U_3 \cup U_4 \cup U_5, \\
%x \in U_4 &\textrm{and}& y \in U_2, \\
%x \in U \setminus U_3 &\textrm{and}& y \in U_3, \\
%{x \in \{v_0, v_2\},} & \textrm{or} & y = v_1, \\
x \in U_1 &\textrm{and}& y \in U \setminus U_1, \\
x \in U_4 &\textrm{and}& y \in U_2, \\
x \in U \setminus U_3 &\textrm{and}& y \in U_3.
\end{array}$
\end{center}

All other edges join vertices within the same $U_i$, and are oriented arbitrarily.
It can be checked that the conditions are self-consistent.
% though not mutually exclusive.
We obtain $D_1$ from $D_2$ by reversing the orientation of the arcs incident to $v_2$.
%We will establish:
\newline

\noindent{\sc Claim 1:} If $S$ is a hull-set in $D_2$, then $I_{D_2}(S) \subseteq I_{D_1}(S-v_1)$.
%\newline

%\noindent{\sc Proof of claim 1:}
Since $S$ is a hull-set for $D_2$, it must contain the extreme vertices $v_0$ and $v_2$. In $D_1$, $v_1$ is on a $v_0 - v_2$ geodesic, and is thus in $I_{D_1}(S-v_1)$. So $S \subseteq I_{D_1}(S-v_1)$.

Consider, therefore, a vertex $w \in I_{D_2}(S) \setminus S$; note that $w \in U$. 
%We need to show that $w \in I_{D_1}(S-v_1)$.
This vertex must be an internal vertex of an $a-b$ geodesic $P$ in $D_2$, for some $a$ and $b$ in $S$. If $a$ and $b$ are both in $U$, then $V(P) \subseteq U$; since the orientation of $G[U]$ is the same in $D_1$ as in $D_2$, $P$ is present in $D_1$. Moreover, the $a-b$ dipaths in $D_1$ are just the $a-b$ dipaths in $D_2$, so $P$ is still a \emph{shortest} $a-b$ dipath. Since $a$ and $b$ are in $S-v_1$, $w \in I_{D_1}(S-v_1)$. 

If $a = v_0$, then $b \not= v_1$ (since the only $v_0-v_1$ geodesic is $\overrightarrow{v_0v_1}$), and clearly $b \not= v_2$, so $b \in U$. Moreover, the $a-b$ dipaths do not use $v_1$ or $v_2$, so $D_1$ contains all the $a-b$ dipaths of $D_2$, and no others; thus $P$ is still an $a-b$ geodesic in $D_1$.
As above, $a$ and $b$ are in $S-v_1$, so $w \in I_{D_1}(S-v_1)$.

If $a = v_2$, then $b$ must be in $N(v_2)$; but then the unique $a-b$ geodesic in $D_2$ is $\overrightarrow{ab}$, with no internal vertices.

If $b= v_1$, then I claim that $P$ must have vertices $awv_1$, with $a \in U_4$ and $w \in U_2$. To see this, note that $a$ cannot be in $N(v_1)$, as otherwise the only $a-v_1$  geodesic is $\overrightarrow{a v_1}$. Moreover, there are no dipaths from $U_3 \cup U_5$ to $v_1$, so $a$ must be in $U_4$. By definition of $U_4$, and by the choice of orientation, there is a (directed) path of length two from $a$ to $v_1$, so every $a-v_1$ geodesic has length two. The internal vertex must be adjacent to $v_1$, but cannot be in $U_1$ (by choice of orientation), so it must be in $U_2$.

Since $a$ is in $U_4$, it is not adjacent to $v_2$; but in $D_1$ there is a directed path $awv_2$, and this is therefore an $a-v_2$ geodesic. Since $a$ and $v_2$ are in $S-v_1$, $w$ is in $I_{D_1}(S-v_1)$.
\newline

\noindent{\sc Claim 2:} If $S$ is a hull-set in $D_2$, then $I^{\ell}_{D_2}(S) \subseteq I^{\ell}_{D_1}(S-v_1)$ for any $\ell \geq 1$.
%\newline

We proceed by induction on $\ell$, the base case $\ell = 1$ following from Claim 1. 
Now for $\ell > 1$, 
\begin{eqnarray*}
I^{\ell}_{D_2}(S) &=& I_{D_2}(I^{\ell-1}_{D_2}(S)) \subseteq
I_{D_1}(I^{\ell-1}_{D_2}(S)-v_1) \subseteq \\
&\subseteq& I_{D_1}(I^{\ell-1}_{D_1}(S-v_1)-v_1) \subseteq I_{D_1}(I^{\ell-1}_{D_1}(S-v_1)) = I^{\ell}_{D_1}(S-v_1).
\end{eqnarray*}
The first containment follows from Claim 1 applied to the hull-set $I^{\ell-1}_{D_2}(S)$, while the second follows from the inductive hypothesis.
\newline

If $S$ is a hull-set for $D_2$, then $I^k_{D_2}(S) = V$, for some $k$. By Claim 2, $I^k_{D_1}(S-v_1) = V$, so $S-v_1$ is a hull-set for $D_1$. 
In particular, $v_1$ is a sink in $D_2$, so it is contained in $S$, and taking $S$ to be a minimum hull-set for $D_2$ we have 
\[ h^{-}(G) \leq h(D_1) \leq |S-v_1| < |S| = h(D_2) \leq h^{+}(G). \] 
If $S$ is a (minimum) geodetic set for $D_2$, then we can take $k=1$, so $S-v_1$ is a geodetic set for $D_1$ and we have $g^{-}(G) < g^{+}(G)$. 
\resnumend
\newline

Since every geodetic set is a hull-set, we have $h(D) \leq g(D)$ for every digraph $D$. For an undirected graph $G$  we therefore have $h^{-}(G) \leq g^{-}(G)$ and $h^{+}(G) \leq g^{+}(G)$, and together with Theorem~\ref{g-h-eqn} this leaves five possibilities:

\begin{eqnarray}
h^{-} = g^{-} &<& h^{+} = g^{+} \label{hg1} \\
h^{-} = g^{-} &<& h^{+} < g^{+} \label{hg3} \\
h^{-} < g^{-} &<& h^{+} = g^{+} \label{hg4} \\
h^{-} < g^{-} &=& h^{+} < g^{+} \label{hg5} \\
h^{-} < h^{+} &<& g^{-} < g^{+} \label{hg6}.
\end{eqnarray}

Chartrand et al. identified many infinite classes of graphs for which~(\ref{hg1}) holds, including trees, cycles and complete bipartite graphs. For complete bipartite graphs $K_{s,t}$ with $s \geq t \geq 2$ \cite[Prop. 3.8]{cz00}, and for transitively orientable graphs with a Hamiltonian path, we have $h^{-}(G) = g^{-}(G) = 2 < n = h^{+}(G) = g^{+}(G)$. 
If $T$ is a tree with $k$ end-vertices, then 
$h^{-}(T) = g^{-}(T) = k < |V(T)| = h^{+}(T) = g^{+}(T)$, while 
$h^{-}(C_{2n+1}) = g^{-}(C_{2n+1}) = 2 < 2n = h^{+}(C_{2n+1}) = g^{+}(C_{2n+1})$. 
We leave the realisability of (\ref{hg3}) -- (\ref{hg6}) as open problems.

\stmtnoboxnum{Problem} 
{Find infinite classes of graphs for which (\ref{hg3}), (\ref{hg4}) or (\ref{hg5}) hold. Are there (infinitely many) graphs for which (\ref{hg6}) holds?
}
\stmtnoboxnumend
\newline

Note that (\ref{hg6}) cannot hold for graphs $G$ for which there is an orientation $\overrightarrow{G}$ such that $g(\overrightarrow{G}) = h(\overrightarrow{G})$. However, there are probably many graphs for which no such orientation exists.

%\stmtnum{Problem} 
%{Is it true that for every graph $G$ there is an orientation $\overrightarrow{G}$ such %that $g(\overrightarrow{G}) = h(\overrightarrow{G})$? This would imply $g^{-}(G) %\leq h^{+}(G)$.
%}
%\stmtnumend
%%\newline

\end{document}

%% file: geo-hull.pstex_t
\begin{picture}(0,0)%
\epsfig{file=geo-hull.pstex}%
\end{picture}%
\setlength{\unitlength}{3947sp}%
\begingroup\makeatletter\ifx\SetFigFont\undefined%
\gdef\SetFigFont#1#2#3#4#5{%
  \reset@font\fontsize{#1}{#2pt}%
  \fontfamily{#3}\fontseries{#4}\fontshape{#5}%
  \selectfont}%
\fi\endgroup%
\begin{picture}(5032,7961)(451,-7473)
\put(2626,-2661){\makebox(0,0)[lb]{\smash{\SetFigFont{12}{14.4}{\rmdefault}{\mddefault}{\updefault}$v_0$}}}
\put(451,-736){\makebox(0,0)[lb]{\smash{\SetFigFont{12}{14.4}{\rmdefault}{\mddefault}{\updefault}$D_1$}}}
\put(451,-5161){\makebox(0,0)[lb]{\smash{\SetFigFont{12}{14.4}{\rmdefault}{\mddefault}{\updefault}$D_2$}}}
\put(5476,-1786){\makebox(0,0)[lb]{\smash{\SetFigFont{12}{14.4}{\rmdefault}{\mddefault}{\updefault}$U_3$}}}
\put(2626,-1786){\makebox(0,0)[lb]{\smash{\SetFigFont{12}{14.4}{\rmdefault}{\mddefault}{\updefault}$U_1$}}}
\put(3676,-1786){\makebox(0,0)[lb]{\smash{\SetFigFont{12}{14.4}{\rmdefault}{\mddefault}{\updefault}$U_2$}}}
\put(4076,-211){\makebox(0,0)[lb]{\smash{\SetFigFont{12}{14.4}{\rmdefault}{\mddefault}{\updefault}$U_4$}}}
\put(2276, 14){\makebox(0,0)[lb]{\smash{\SetFigFont{12}{14.4}{\rmdefault}{\mddefault}{\updefault}$U_5$}}}
\put(4626,-2661){\makebox(0,0)[lb]{\smash{\SetFigFont{12}{14.4}{\rmdefault}{\mddefault}{\updefault}$v_2$}}}
\put(2626,-7086){\makebox(0,0)[lb]{\smash{\SetFigFont{12}{14.4}{\rmdefault}{\mddefault}{\updefault}$v_0$}}}
\put(4626,-7086){\makebox(0,0)[lb]{\smash{\SetFigFont{12}{14.4}{\rmdefault}{\mddefault}{\updefault}$v_2$}}}
\put(3626,-7086){\makebox(0,0)[lb]{\smash{\SetFigFont{12}{14.4}{\rmdefault}{\mddefault}{\updefault}$v_1$}}}
\put(3626,-2661){\makebox(0,0)[lb]{\smash{\SetFigFont{12}{14.4}{\rmdefault}{\mddefault}{\updefault}$v_1$}}}
\put(5476,-6211){\makebox(0,0)[lb]{\smash{\SetFigFont{12}{14.4}{\rmdefault}{\mddefault}{\updefault}$U_3$}}}
\put(2626,-6211){\makebox(0,0)[lb]{\smash{\SetFigFont{12}{14.4}{\rmdefault}{\mddefault}{\updefault}$U_1$}}}
\put(3676,-6211){\makebox(0,0)[lb]{\smash{\SetFigFont{12}{14.4}{\rmdefault}{\mddefault}{\updefault}$U_2$}}}
\put(4076,-4636){\makebox(0,0)[lb]{\smash{\SetFigFont{12}{14.4}{\rmdefault}{\mddefault}{\updefault}$U_4$}}}
\put(2276,-4411){\makebox(0,0)[lb]{\smash{\SetFigFont{12}{14.4}{\rmdefault}{\mddefault}{\updefault}$U_5$}}}
\end{picture}